\documentclass{amsart}
\usepackage{amsmath, amsthm, amscd, amsfonts}

    \newtheorem{prop}{Proposition}[section]
    \newtheorem{cor}[prop]{Corollary}
    \newtheorem{lem}[prop]{Lemma}
    \newtheorem{rem}[prop]{Remark}
    \newtheorem{defn}[prop]{Definition}
    \newtheorem{thm}[prop]{Theorem}

  \renewcommand{\t}[1]{\mbox{$\tilde{#1}$}}

  \newcommand{\jbst}{$JB^*$-triple}
 \newcommand{\jbwst}{$JBW^*$-triple}

 \newcommand{\csa}{$C^*$-algebra}

  \newcommand{\pf}{{\it Proof}.}

   \newcommand{\tp}[3]{\{#1#2#3\}}
\newcommand{\tprime}[3]{\{#1#2#3\}^{\prime}}
   \newcommand{\tpc}[3]{\{#1,#2,#3\}}
\newcommand{\ter}[3]{[#1#2#3]}
\newcommand{\terc}[3]{[#1,#2,#3]}
\newcommand{\tep}[3]{(#1#2#3)}
\newcommand{\tea}[3]{\langle #1#2#3\rangle}
\newcommand{\tepc}[3]{(#1,#2,#3)}
\newcommand{\teac}[3]{\langle #1,#2,#3\rangle}

\newcommand{\bfa}{\mbox{$\bf a$}}
\newcommand{\bfab}{\mbox{$\overline{\bf a}$}}
\newcommand{\bfb}{\mbox{$\bf b$}}
\newcommand{\bfbb}{\mbox{$\overline{\bf b}$}}

\begin{document}
\title[Operator space characterizations of C*-algebras and 
ternary rings]
{Operator space characterizations of C*-algebras
and  ternary rings \footnote{This work was
supported in part by NSF grant DMS-0101153}}
\author{Matthew Neal \and Bernard Russo}
\address{
Denison University,
Granville, Ohio 43023\newline
\indent        University of California,
 Irvine, California 92697-3875}
\email{nealm@denison.edu brusso@math.uci.edu}
\begin{abstract}
We prove that an operator space is completely isometric to a ternary ring
of operators if and only if the open unit balls of all of its matrix spaces
are bounded symmetric domains. From this we obtain an operator
space characterization of C*-algebras.
\end{abstract}
\maketitle


\begin{center}
{\bf Introduction}
\end{center}
In the category of operator spaces, that is, 
subspaces of the bounded linear operators $B(H)$
on a complex Hilbert space $H$ together with the
induced matricial operator norm structure, objects are
equivalent if they are completely isometric, {\it i.e.} if there is a linear
isomorphism between the spaces which preserves this
matricial norm structure. Since operator algebras, that is,
subalgebras of $B(H)$, are  motivating
examples for much of operator space theory, it is natural
to ask if one can characterize which  operator spaces are
operator algebras. One satisfying answer was given by Blecher,
Ruan and Sinclair in \cite{BL90}, where it was shown that among operator
spaces $A$ with a (unital but not necessarily associative)
Banach algebra product, those which are
completely isometric to operator algebras are precisely the ones whose
multiplication is completely contractive with respect to the
Haagerup norm on $A\otimes A$
(For a completely bounded version of this 
result, see \cite{Blecher95}).  

A natural object to characterize in this context are
the so called {\it ternary rings of operators} (TRO's).
These are subspaces of $B(H)$ which
are closed under the ternary
product $xy^{\ast}z$. This class includes
C*-algebras. 
TRO's, like C*-algebras, carry a natural operator space structure.
In fact, every TRO is (completely)
 isometric to a corner
$pA(1-p)$ of a C*-algebra $A$. TRO's are important because, as shown by
Ruan \cite{Ruan89}, the injectives in the category of
operator spaces are TRO's (corners of injective C*-algebras) and not, in
general, operator algebras (For the dual version of this result see
\cite{EffOzaRua}).
Injective envelopes of operator systems and 
of operator spaces (\cite{Ham79},\cite{Ruan89})
 have proven to be important tools,
see for example \cite{BlePau01}. The
characterization of TRO's among operator spaces is the subject of this paper
(See Theorem~\ref{thm:first}).

Closely related to TRO's are the so called JC*-triples, norm closed
subspaces of $B(H)$
which are closed under the triple product
$(xy^{\ast}z + zy^{\ast}x)/2$.
 These generalize the class of TRO's and
have the property, as shown by Harris in \cite{Harris73},
that isometries coincide with algebraic isomorphisms. It is not hard to
see this implies that the algebraic isomorphisms in the
class of TRO's are complete isometries, since for each TRO
$A$, $M_{n}(A)$ is a JC*-triple (For the converse of this, see
\cite[Proposition 2.1]{Ham99}).
As a consequence, if an operator space $X$ is completely isometric
to a TRO, then the induced ternary product on $X$ is unique, {\it i.e.},
independent of the TRO.

Building on the pioneering work of Arveson
 (\cite{Arveson69},\cite{Arveson72}) on noncommutative analogs of 
the Choquet and Shilov boundaries,
Hamana (see
\cite{Ham99}) proved that every operator space $A$ has a unique enveloping
TRO ${\mathcal T}(A)$ which is an invariant of complete
isometry and has the property that for
any TRO $B$ generated by a realization of $A$,
there exists a homomorphism of $B$ onto ${\mathcal T}(A)$. 
The space ${\mathcal T}(A)$
is also called the {\it Hilbert C*-envelope}
of $A$. The work in \cite{Blecher01} suggests
that the Hilbert C*-envelope
is an appropriate noncommutative generalization to operator spaces
of the classical theory
of Shilov boundary of function spaces.

It is also true that
a commutative TRO ($xy^*z=zy^*x$) is an
 associative JC*-triple and hence by \cite[Theorem 2]{FriRus83},
 is isometric (actually completely
isometric) to a complex $C_{hom}$-space, that is, the
space of weak*-continuous functions on 
the set of extreme points of the unit ball of the dual
of a Banach space which are homogeneous with respect to
the natural action of the circle group, see \cite{FriRus83}.
Hence, if one views operator spaces as
noncommutative Banach spaces, and C*-algebras as noncommutative $C(\Omega)$'s, then
TRO's and JC*-triples
can be viewed as noncommutative $C_{hom}$-spaces.

As noted above, injective operator spaces, {\it i.e.}, those
which are the range of a {\it completely}
contractive projection on some $B(H)$, are completely isometrically
TRO's; the so called {\it mixed injective} operator spaces, those which are
the range of a contrative projection on some
$B(H)$, are isometrically JC*-triples. The operator space 
classification of mixed injectives
was begun by the authors in \cite{NR} and is ongoing.

Relevant to this paper is another property shared by all JC*-triples (and
hence all TRO's). For any Banach space $X$, we denote by $X_0$ its
open unit ball: $\{x\in X:\|x\|<1\}$. 
The open unit ball of every JC*-triple is a {\it bounded
symmetric domain}. This is equivalent to saying that it has a
transitive group of biholomorphic automorphisms. It was shown by Koecher
in finite dimensions (see \cite{Loos77}) and Kaup
\cite{Kaup83} in the general case that this is a defining property for the
slightly larger class of JB*-triples. The only
illustrative basic examples of  JB*-triples which are not  JC*-triples
are the
space $H_{3}({\mathcal O})$ of 3 x 3 hermitian matrices
over the octonians and a certain subtriple of $H_{3}({\mathcal O})$. These
are called  {\it exceptional} triples, and they
cannot be represented as a JC*-triple. This holomorphic characterization
has been useful as it gives an elegant proof, due to Kaup
\cite{Kaup84}, that the range of a contractive projection on a JB*-triple
is isometric to another JB*-triple. The same statement
holds for JC*-triples, as proven earlier 
by Friedman and Russo in \cite{FriRus85}.
Youngson proved in
\cite{Youngson83} that the range of a completely contractive projection on
a C*-algebra is completely isometric to a TRO. These results, as
well as those of \cite{AraFri78} and \cite{EffSto79}, are rooted
in the fundamental result of Choi-Effros \cite{ChoEff77} for
completely positive projections on C*-algebras and the classical
result (\cite{LinWul69},\cite[Theorem 5]{FriRus82})
 that the range of a contractive projection on $C(\Omega)$ is isometric to
a $C_\sigma$-space, hence a $C_{hom}$-space.

Motivated by this characterization for JB*-triples, we will give a
holomorphic characterization of TRO's
up to complete isometry. We will prove in Theorem~\ref{thm:first}
that an operator space $A$ is
completely isometric to a TRO if and only if the open unit
balls $M_{n}(A)_{0}$ are bounded symmetric domains for all $n \geq 2$. As a
consequence, we obtain in Theorem~\ref{thm:4.5} a holomorphic operator
space characterization of C*-algebras as well. It should be mentioned that
Upmeier (for the category of Banach spaces) in \cite{Upmeier84} and El Amin-Campoy-Palacios 
(for the category of Banach algebras) in \cite{AmiCamPal01},
gave  different but still
holomorphic characterizations of C*-algebras up to isometry.
We note in passing that injective operator spaces satisfy
the  hypothesis of Theorem~\ref{thm:first}, 
so we obtain that they are (completely isometrically)
TRO's based on deep results about JB*-triples rather than the
deep result of Choi-Effros (See Corollaries~\ref{cor:4.4} and ~\ref{cor:4.6}).

We now describe the organization of this paper. Section~\ref{sect:prel} 
contains the necessary background and some preliminary results on
 contractive projections. In section~\ref{sect:additivity}, three 
auxiliary ternary products are introduced and are shown to yield the original
JB*-triple product upon symmetrization. Section~\ref{sect:equality} is
devoted to proving that these three ternary products all coincide. 
Section~\ref{sect:main} contains the statement and proof of the
main result and its consequences.

The authors wish to thank Profesors Zhong-Jin
Ruan and David Blecher for their advice and encouragement at the beginning 
stages of this work.

\section{Preliminaries}\label{sect:prel}

An {\it \bf{operator space}} will be defined as a normed space $A$
together with a linearly isometric representation as a subspace
of some $B(H)$. This gives $A$ a family of operator norms $\| \cdot
\|_{n}$ on $M_{n}(A) \subset B(H^{n})$. As proved in \cite{Ruan88},
an operator space can also be defined abstractly as a normed space $A$
having a norm on $M_{n}(A)$ ($n\ge 2$) satisfying certain
properties. Each such family of norms is regarded as a
``quantization'' of the underlying Banach space.
These properties give rise to an isometric representation of
the operator space as a subspace of B(H) where the
natural amplification maps preserve the matricial norm structure. This is
analagous to (and generalizes) the way an abstract Banach
space $B$ can be isometrically embedded as a subspace of $C(\Omega)$. The
resulting operator space structure in this case is called
$MIN(B)$ and is seen as a commutative quantization of $B$.

Two operator spaces $A$ and $B$ are {\it \bf{n-isometric}} if there exists
an isometry $\phi$ from $A$ onto $B$
such that the amplification mapping $\phi_{n}:M_{n}(A) \rightarrow
M_{n}(B)$ defined by $\phi([a_{ij}])=[\phi(a_{ij})]$ is an
isometry. $A$ and $B$ are {\it \bf{completely isometric}} if there exists
a mapping $\phi$ from $A$ onto $B$ which is an
n-isometry for all n. For other basic results about operator spaces, see
\cite{EffRua00}.

The following definition is
a Hilbert space-free generalization of the TRO's mentioned in the
introduction. 
\begin{defn}[Zettl \cite{Zettl83}] A  
{\it \bf{C*-ternary ring}} is a Banach space $A$ with
ternary product $[x,y,z]:A\times A\times A\rightarrow A$ 
which is linear in the outer variables,
conjugate linear in the middle variable, is {\bf associative}:
\[
[ab[cde]]=[a[dcb]e]=[ab[cde]],
\]
and satisfies $\|[xyz]\|\le\|x\|\|y\|\|z\|$ and  $\|[xxx]\|=\|x\|^{3}$.
\end{defn}

A TRO is a C*-ternary ring under any of the  products $\ter{x}{y}{z}_\lambda=
\lambda xy^*z$, for any complex number $\lambda$ with $|\lambda|= 1$.

A linear map $\varphi$ between C*-ternary rings is a {\bf homomorphism} if
$\varphi(\ter{x}{y}{z})=\terc{\varphi(x)}{\varphi(y)}{\varphi(z)}$ and an
{\bf anti-homomorphism} if
$\varphi(\ter{x}{y}{z})=-\terc{\varphi(x)}{\varphi(y)}{\varphi(z)}$.

The following is a Gelfand-Naimark representation
theorem for C*-ternary rings.  

\begin{thm}[\cite{Zettl83}]\label{thm:1} 
 For any C*-ternary ring $A$,
$A=A_{1} \oplus A_{-1}$,
where 
$A_1$ and $A_{-1}$ are sub-C*-ternary rings,
$A_{1}$ is isometrically isomorphic to a TRO $B_{1}$ and $A_{-1}$ is
isometrically anti-isomorphic  to a TRO
$B_{-1}$.  
\end{thm}
\medskip

It follows that $A_{-1}=0$ if and only if $A$ is ternary isomorphic
to a TRO.  In
 Theorem~\ref{thm:first}, we shall show that under suitable
assumptions on an operator space $A$, it becomes a C*-ternary ring
with $A_{-1}=0$ and the above ternary isomorphism is a complete
isometry from $A$  with its original operator
space structure to a TRO with its natural operator space
structure.

An immediate consequence of our proof of  Theorem~\ref{thm:first} is
an answer to a question posed by Zettl \cite[p.136]{Zettl83}:
for a C*-ternary ring $A$, $A_{-1}=0$ if and only if $A$ is a
JB*-triple (see the next definition) under the triple product 
\[
\tp{a}{b}{c}=\frac{1}{2}(\ter{a}{b}{c}+\ter{c}{b}{a}). 
\]

\medskip

The following definition generalizes 
the JC*-triples defined in the introduction.

\begin{defn}[\cite{Kaup83}]\label{defn:1.2} 
A {\it {\bf JB*-triple}} is a Banach space $A$ with a product
$D(x,y)z = \{ x \,\ y \,\ z \}$ which is linear in
the outer variables, conjugate linear in the middle variable, is {\bf
commutative}: $\{ x \,\ y \,\ z \} = \{ z \,\ y \,\ x \}$,
satisfies an  associativity condition:
\begin{equation}\label{eq:main}
[D(x,y),D(a,b)] = D(\{ x \,\ y \,\ a \},b) - D(a,\{ b \,\ x \,\ y \})
\end{equation}
and has the topological properties that {\bf (1)} $\|D(x,x)\|=\|x\|^{2}$
{\bf (2)} $D(x,x)$ is {\bf hermitian} (in the sense that
$\|e^{itD(x,x)}\|=1$) and has positive spectrum in the Banach algebra
$B(A)$.  We abbreviate $D(x,x)$ to $D(x)$.
\end{defn}

As noted in the introduction, JC*-triples, (and hence TRO's and
C*-algebras) are examples of JB*-triples. Other examples include
any Hilbert space, and the spaces of symmetric and anti-symmetric elements
of $B(H)$ under a transpose map defined by a conjugation.

If one ignores the norm and the topological properties in 
Definition~\ref{defn:1.2}, the algebraic structure which results,
called a {\bf Jordan triple system}, or {\bf Jordan pair},
has a life of its own, \cite{Loos77}.  Note that (\ref{eq:main}) can
be written as

\begin{equation}\label{eq:JP15}
\tpc{x}{y}{\tp{a}{b}{z}}-\tpc{a}{b}{\tp{x}{y}{z}}=\tpc{\tp{x}{y}{a}}{b}{z}-
\tpc{a}{\tp{y}{x}{b}}{z}.
\end{equation}

For easy reference we record here
two identities for Jordan triple systems which can be derived from 
(\ref{eq:main}) (\cite[JP8,JP16]{Loos77}).

\begin{equation}\label{eq:JP8}
2D(x,\tp{y}{x}{z})=D(\tp{x}{y}{x},z)+D(\tp{x}{z}{x},y)
\end{equation}

\begin{equation}\label{eq:JP16}
\tpc{\tp{x}{y}{a}}{b}{z}-\tpc{a}{\tp{y}{x}{b}}{z}=\tpc{x}{\tp{b}{a}{y}}{z}
-\tpc{\tp{a}{b}{x}}{y}{z}
\end{equation}

We will now list some facts about JB*-triples that are relevant to our
paper. A survey of the basic theory can be found in
\cite{Russo94}. As proved by Kaup \cite{Kaup83}, JB*-triples are in 1-1
isometric correspondence with Banach spaces whose open unit
ball is a {\it {\bf bounded symmetric domain}}. The triple product here
arises from the Lie algebra of the group of biholomorphic
automorphisms. This Lie algebra is the space of complete vector fields on
the open unit ball and consists of certain polynomials of degree 
at most 2. 
The quadratic term in each of these polynomials is determined by the
constant term. For a bounded symmetric domain, the constant terms
which occur exhaust $A$. Thus, linearizing the quadratic term for
every element $a\in A$ leads to
a triple product on $A$.

It is this correspondence which motivates the study of the more general
JB*-triples.
Indeed, the proofs of two important facts follow naturally from the
holomorphic point of view \cite{Kaup84}. Firstly, the isometries
between JB*-triples are precisely the algebraic isomorphisms. From this
follows the important fact, used several times in this paper,
 that, unlike the case for binary products, 
 the triple product of a JB*-triple is unique. Secondly,
the range of a contractive projection $P$ on a JB*-triple $Z$ is
isometric to a JB*-triple.  More precisely, $P(Z)$ is a JB*-triple
under the norm and linear operations it inherits from $Z$ 
and the triple product
$\tp{x}{y}{z}_{P(Z)}:=P(\tp{x}{y}{z}_Z)$, for $x,y,z\in P(Z)$.

In the context of JC*-triples, these facts were
proven by functional analytic methods in \cite{Harris73}
and \cite{FriRus85} respectively. These facts show that JB*-triples are a
natural category in which to study isometries and
contractive projections. Recently, in \cite{CNR} the authors with C-H. Chu
have shown that w*-continuous contractive projections on
dual JB*-triples (called JBW*-triples) preserve the Jordan
triple generalization
of the Murray-von-Neumann type decomposition established
in \cite{Horn87} and \cite{H1}. Two other properties of contractive
projections were used in that work and will be needed in the present
paper. They consist of two conditional expectation formulas 
for contractive projections on JC*-triples (\cite[Corollary 1]{FriRus84})
\begin{equation}\label{eq:CE}
P\tpc{Px}{Py}{Pz}=P\tpc{Px}{Py}{z}=P\tpc{Px}{y}{Pz};
\end{equation}
and the
fact that the range of a {\it bi}contractive projection on a JC*-triple is
a subtriple \cite[Proposition 1]{FriRus84}.

 Let $A$ be a JB*-triple. For any $a \in A$, there is a triple functional
calculus, that is,
a triple isomorphism of the closed subtriple
 $C(a)$ generated by $a$ onto the commutative
C*-algebra $C_{0}(  \mbox{Sp}\, D(a,a)\cup\{0\} )$ of continuous functions
vanishing at zero, with the triple product $f\overline{g}h$. Any
JBW*-triple (defined above) has the propertly that it is the
norm closure of the linear span of its {\it \bf {tripotents}}, that is,
elements $e$ with $e=\tp{e}{e}{e}$. A {\it \bf{unitary}} tripotent is
a tripotent $v$ such that $D(v,v)=Id$. For a C*-algebra, tripotents are
the partial isometries and for unital C*-algebras,
unitary tripotents are precisely
the unitaries. For tripotents $u$ and $v$, algebraic orthogonality,
 i.e. $D(u,v)=0$,
coincides with Banach space othogonality:
$\|u \pm v\|=1$. For $a$ and $b$ in $A$, we will denote the
property $D(a,b)=0$ by $a \perp b$.

As proved in \cite{D}, the second dual $A^{\ast\ast}$ of a JB*-triple
$A$ is a JBW*-triple containing $A$ as a subtriple. Multiplication in a
JBW*-triple is  norm
continuous and, as proved in \cite{BarTim86}, separately w*-continuous. 

We close this section of preliminaries with an elementary proposition
showing that certain concrete projections are contractive.

\begin{prop}\label{prop:1.3}
Let $A$ be an operator space in $B(H)$.
\begin{description}
\item[(a)] Define a projection $P$ on $M_2(A)$ by
\[
P(\left[\begin{array}{cc} a&b\\
c&d\end{array}\right])=\frac{1}{2}\left[\begin{array}{cc} a+b&a+b\\
0&0\end{array}\right].
\]
Then $\|P\|\le 1$. Moreover, the restriction of $P$ to
$\{\left[\begin{array}{cc} a&b\\
0&0\end{array}\right]:a,b\in A\}$
is bicontractive.
\item[(b)]
Let $P_{11}:M_2(A)\rightarrow M_2(A)$ be the map  $\left[\begin{array}{cc}
a_{11}&a_{12}\\
a_{21}&a_{22}\end{array}\right]\mapsto \left[\begin{array}{cc} a_{11}&0\\
0&0\end{array}\right]$, and similarly for $P_{12},P_{21},P_{22}$. Then
$P_{ij}$ is
contractive and $P_{11}+P_{21}$, $P_{11}+P_{12}$, and $P_{11}+P_{22}$ are
bicontractive. More generally, the $P_{ij}:M_n(A)\rightarrow M_n(A)$ are
contractive and for any subset $S\subset \{1,2,\ldots,n\}$,
$\sum_{i\in S}\sum_{j=1}^n P_{ij}$ and $\sum_{j\in S}\sum_{i=1}^n
 P_{ij}$ are bicontractive.

\item[(c)]
The projections $P:M_2(A)\rightarrow M_2(A)$ and $Q:P(M_2(A))\rightarrow
P(M_2(A))$ defined by
\[
P(\left[\begin{array}{cc} a&b\\
c&d\end{array}\right])=\frac{1}{2}\left[\begin{array}{cc} a+d&b+c\\
b+c&a+d\end{array}\right]
\]
and
\[
Q(\left[\begin{array}{cc} a&b\\
b&a\end{array}\right])=\frac{1}{2}\left[\begin{array}{cc} a+b&a+b\\
a+b&a+b\end{array}\right]
\]
are bicontractive.
\end{description}
\end{prop}
\noindent\pf\
We omit the proofs of (a) and (b). To prove (c), since for example
$I-P=(I-(2P-I))/2$ and $P=(I+(2P-I))/2$, it
 suffices to show that $2P-I$ and $2Q-I$ are contractive. But
\[
(2P-I)(\left[\begin{array}{cc} a&b\\ c&d\end{array}\right])=
\left[\begin{array}{cc} d&c\\ b&a\end{array}\right]=
\left[\begin{array}{cc} 0&1\\ 1&0\end{array}\right]
\left[\begin{array}{cc} a&b\\ c&d\end{array}\right]
\left[\begin{array}{cc} 0&1\\ 1&0\end{array}\right],
\]
and
\[
(2Q-I)(\left[\begin{array}{cc} a&b\\ b&a\end{array}\right])=
\left[\begin{array}{cc} b&a\\ a&b\end{array}\right]=
\left[\begin{array}{cc} 0&1\\ 1&0\end{array}\right]
\left[\begin{array}{cc} a&b\\ b&a\end{array}\right].\qed
\]

\section{Additivity of the ternary products}\label{sect:additivity}

Throughout this section, $A\subset B(H)$ will 
be an operator space such that the open unit
ball $M_2(A)_0$ is a bounded symmetric domain. 
Let
$\tpc{\cdot}
{\cdot}{\cdot}_{M_2(A)}$ denote the associated \jbst\ product on
$M_2(A)$. Note that although $M_2(A)$ inherits the norm and linear
structure of $M_2(B(H))=B(H\oplus H)$, its triple product 
$\tp{\cdot}
{\cdot}{\cdot}_{M_2(A)}$ in general differs from the concrete triple product
$(XY^*Z+ZY^*X)/2$ of $B(H\oplus H)$. In fact, the results of this
section would become trivial if these two triple products were the same.

By properties of 
contractive projections and the uniqueness of the triple product,
$A$, being linearly isometric to $P_{ij}(M_2(A))$ becomes a JB*-triple
whose triple product $\tp{x}{y}{z}_A$ is given, for example, by
\[ 
\left[\begin{array}{cc}\tp{x}{y}{z}_A&0\\
0&0\end{array}\right]
=P_{11}\left(\tp{\left[\begin{array}{cc} x&0\\ 0&0\end{array}\right]}
{\left[\begin{array}{cc} y&0\\ 0&0\end{array}\right]}
{\left[\begin{array}{cc} z&0\\ 0&0\end{array}\right]}_{M_2(A)}\right),
\]
and similarly using the other $P_{ij}$.
Usually we shall just
use the notation $\tp{\cdot}{\cdot}{\cdot}$ for either of the triple
products $\tp{x}{y}{z}_A$ and $\tpc{\cdot}
{\cdot}{\cdot}_{M_2(A)}$. Lemma~\ref{lem:2.8} shows that the projection
$P_{11}$ could be removed in this definition.

We assume $A$ is as above and proceed to 
define (in Definition~\ref{defn:2.7},
 three auxiliary ternary products, denoted
$\terc{\cdot}{\cdot}{\cdot}$, $\tepc{\cdot}{\cdot}{\cdot}$, and 
$\teac{\cdot}{\cdot}{\cdot}$ and
show their relation to $\tpc{\cdot}{\cdot}{\cdot}$. 
We begin with a
sequence of lemmas which establish some properties of the terms in the
following identity, where $a,b,c\in A$.

\begin{eqnarray}\label{eq:*}
\tp{\left[\begin{array}{cc} a&a\\
0&0\end{array}\right]}{\left[\begin{array}{cc} 0&b\\
0&0\end{array}\right]}{\left[\begin{array}{cc}
c&c\\ 0&0\end{array}\right]}&=&
\tp{\left[\begin{array}{cc} a&0\\
0&0\end{array}\right]}{\left[\begin{array}{cc} 0&b\\
0&0\end{array}\right]}{\left[\begin{array}{cc}
c&0\\ 0&0\end{array}\right]}\\\nonumber
&+&
\tp{\left[\begin{array}{cc} a&0\\
0&0\end{array}\right]}{\left[\begin{array}{cc} 0&b\\
0&0\end{array}\right]}{\left[\begin{array}{cc}
0&c\\ 0&0\end{array}\right]}\\\nonumber
&+&
\tp{\left[\begin{array}{cc} 0&a\\
0&0\end{array}\right]}{\left[\begin{array}{cc} 0&b\\
0&0\end{array}\right]}{\left[\begin{array}{cc}
0&c\\ 0&0\end{array}\right]}\\\nonumber
&+&
\tp{\left[\begin{array}{cc} 0&a\\
0&0\end{array}\right]}{\left[\begin{array}{cc} 0&b\\
0&0\end{array}\right]}{\left[\begin{array}{cc}
c&0\\ 0&0\end{array}\right]}.
\end{eqnarray}

It will be shown in Lemma~\ref{lem:2.3} that the left side of (\ref{eq:*})
has the form
\[
\left[\begin{array}{cc} x&y\\ z&w\end{array}\right],
\]
where $(x+y)/2=\tp{a}{b}{c}$.
In Lemmas~\ref{lem:2.6}-\ref{lem:2.8}, each term on the right side of
(\ref{eq:*}) will be analyzed.

\begin{rem}\label{rem:2.1}
The space
\[
\t{A}=\{\t{a}=\left[\begin{array}{cc} a&a\\ 0&0\end{array}\right]:a\in A\}
\]
with the triple product
\begin{equation}\label{eq:999}
\tp{\t{a}}{\t{b}}{\t{c}}_{\t{A}}:
=\left[\begin{array}{cc}
2\tp{a}{b}{c}&2\tp{a}{b}{c}\\ 0&0\end{array}\right]
\end{equation}
and the norm of $M_2(A)$, is a \jbst. 
\end{rem}

Note that by 
Proposition~\ref{prop:1.3}(a), $\tilde{A}$ is a subtriple of
$M_2(A)$, but we do not know {\it a priori} that its triple
product is given by (\ref{eq:999}).

\noindent\pf\   
The proposed triple product, which we denote by
$\tp{\t{a}}{\t{b}}{\t{c}}$,  is obviously linear and symmetric in $\t{a}$
and $\t{c}$, and conjugate linear
in $\t{b}$.  Since, for example,
\[
\tp{\t{a}}{\t{b}}{\tp{\t{c}}{\t{d}}{\t{e}}}=\left[\begin{array}{cc}
2\tp{a}{b}{\tp{c}{d}{e}}&2\tp{a}{b}{\tp{c}{d}{e}}\\
0&0\end{array}\right],
\]
the main identity (\ref{eq:JP15}) is satisfied.

From $\|\left[\begin{array}{cc} a&a\\
0&0\end{array}\right]\|=\sqrt{2}\|a\|$ one obtains
 $\|\tp{\t{a}}{\t{a}}{\t{a}}\|=\|\t{a}\|^3$,
$\|\tp{\t{a}}{\t{b}}{\t{c}}\|\le \|\t{a}\||\|\t{b}\|\|\t{c}\|$ and
hence $\|D(\t{a})\|=\|\t{a}\|^2$.

Since $e^{itD(\t{x})}\t{y}=(e^{2itD(x)}y)^{\t{}}$,
$
\|e^{itD(\t{x})}\t{y}\|=\sqrt{2}\|e^{2itD(x)}y\|=
\sqrt{2}\|y\|=\|\t{y}\|,
$
so $D(\t{x})$ is hermitian.

Finally, for $\lambda<0$, the inverse of $\lambda-D(\t{x})$ 
is given by
\[
\t{y}\mapsto \left[\begin{array}{cc}
(\lambda-2D(x))^{-1}y&(\lambda-2D(x))^{-1}y\\ 0&0\end{array}\right].
\]
Hence, $Sp_{B(\tilde{A})}(D(\t{x}))\subset [0,\infty)$.\qed

\begin{lem}\label{lem:2.3} For $a,b,c\in A$, there exist $x,y,z,w\in A$
such that
\[
\tp{\left[\begin{array}{cc} a&a\\
0&0\end{array}\right]}{\left[\begin{array}{cc} 0&b\\
0&0\end{array}\right]}{\left[\begin{array}{cc} c&c\\
0&0\end{array}\right]}
=\left[\begin{array}{cc} x&y\\ z&w\end{array}\right],
\]
and $(x+y)/2=\tp{a}{b}{c}$.
\end{lem}
\noindent\pf\
Consider the projection $P$ defined in Proposition~\ref{prop:1.3}(a). By
(\ref{eq:CE}),
\[
P(\tp{\left[\begin{array}{cc} a&a\\
0&0\end{array}\right]}{\left[\begin{array}{cc} 0&b\\
0&0\end{array}\right]}{\left[\begin{array}{cc} c&c\\
0&0\end{array}\right]})=
P(\tp{\left[\begin{array}{cc} a&a\\
0&0\end{array}\right]}{\left[\begin{array}{cc} b/2&b/2\\
0&0\end{array}\right]}{\left[\begin{array}{cc} c&c\\
0&0\end{array}\right]}).
\]

By Remark~\ref{rem:2.1} and the uniqueness of the triple product in a
\jbst,
\[
P(\tp{\left[\begin{array}{cc} a&a\\
0&0\end{array}\right]}{\left[\begin{array}{cc} b&b\\
0&0\end{array}\right]}{\left[\begin{array}{cc} c&c\\
0&0\end{array}\right]})
=2\left[\begin{array}{cc} \tp{a}{b}{c}&\tp{a}{b}{c}\\
0&0\end{array}\right].
\]
Thus, if
\[
\tp{\left[\begin{array}{cc} a&a\\
0&0\end{array}\right]}{\left[\begin{array}{cc} 0&b\\
0&0\end{array}\right]}{\left[\begin{array}{cc} c&c\\
0&0\end{array}\right]}=\left[\begin{array}{cc} x&y\\
z&w\end{array}\right],
\]
then
\[
\left[\begin{array}{cc} \tp{a}{b}{c}&\tp{a}{b}{c}\\
0&0\end{array}\right]=\left[\begin{array}{cc} (x+y)/2&(x+y)/2\\
0&0\end{array}\right]. \qed
\]

It will be shown below in the proof of Lemma~\ref{1stprod}
that $x=y=\tp{a}{b}{c}$ and that each $z=w=0$.

\begin{lem}\label{lem:orth}
For each $a,b\in A$,
\[
\left[\begin{array}{cc} a&0\\
0&0\end{array}\right]\perp\left[\begin{array}{cc} 0&0\\
0&b\end{array}\right]
\mbox{ and }
\left[\begin{array}{cc} 0&0\\
a&0\end{array}\right]\perp\left[\begin{array}{cc} 0&b\\
0&0\end{array}\right]
\]
\end{lem}
\noindent\pf\
Suppose first that $a=\sum\lambda_i u_i$ where $\lambda_i>0$ and the $u_i$
are tripotents in $A$, and
similarly for $b=\sum\mu_j v_j$. Because the image of a bicontractive
projection is a subtriple (\cite[Proposition 1]{FriRus84}),
$U_i:=\left[\begin{array}{cc} u_i&0\\ 0&0\end{array}\right]$ and
$V_j:=\left[\begin{array}{cc} 0&0\\ 0&v_j\end{array}\right]$ are
tripotents, and since they are orthogonal in $B(H\oplus H)$,
$\|U_i\pm V_j\|=1$. Hence $D(U_i, V_j)=0$ in (the abstract triple
product of) $M_2(A)$ and so
for all $x,y,z,w\in A$,
\[
\tp{\left[\begin{array}{cc} a&0\\
0&0\end{array}\right]}{\left[\begin{array}{cc} 0&0\\
0&b\end{array}\right]}{\left[\begin{array}{cc} x&y\\
z&w\end{array}\right]}=
\sum_{i,j}\lambda_i\mu_j
\tp{\left[\begin{array}{cc} u_i&0\\
0&0\end{array}\right]}{\left[\begin{array}{cc} 0&0\\
0&v_j\end{array}\right]}{\left[\begin{array}{cc}x&y\\
z&w\end{array}\right]}=0.
\]
For the general case, note that, by \cite[3.2.1]{EffRua00}, 
there is an operator space
structure on the dual of any operator space $A$
such that the canonical inclusion of $A$ into $A^{**}$
is a complete isometry. 
Moreover, by \cite[Theorem 2.5]{Blecher92}
the norm structure on $M_{n}(A^{\ast\ast})$
coincides with that obtained from the identification $M_{n}(A^{\ast\ast})
= M_{n}(A)^{\ast\ast}$. Hence, for all $n$, $M_{n}(A^{\ast\ast})$ is
a JBW*-triple containing $M_{n}(A)$ as subtriple. Since each element of
$A$ can be approximated in norm by finite linear
combinations of tripotents in $A^{\ast\ast}$, the first statement in the
lemma follows from the norm continuity of the triple product.

Since interchanging rows is an isometry, hence an isomorphism, the second
statement follows.
\qed

\begin{lem}\label{lem:2.6} Let $a,b,c\in A$. Then
\[ 
\tp{\left[\begin{array}{cc} a&0\\
0&0\end{array}\right]}{\left[\begin{array}{cc} 0&b\\
0&0\end{array}\right]}{\left[\begin{array}{cc} c&0\\
0&0\end{array}\right]}=0,
\] 
\begin{equation}\label{eq:7bis}
\tp{\left[\begin{array}{cc} a&0\\
0&0\end{array}\right]}{\left[\begin{array}{cc} 0&0\\
b&0\end{array}\right]}{\left[\begin{array}{cc} c&0\\
0&0\end{array}\right]}=0,
\end{equation}.
\end{lem}
\noindent\pf\
To prove
the first statement, let
 $X$ denote $\tp{\left[\begin{array}{cc} a&0\\
0&0\end{array}\right]}{\left[\begin{array}{cc} 0&b\\
0&0\end{array}\right]}{\left[\begin{array}{cc} c&0\\
0&0\end{array}\right]}$.
By (\ref{eq:CE}),
\[
P_{11}(X)=P_{11}(\tp{\left[\begin{array}{cc} a&0\\
0&0\end{array}\right]}{\left[\begin{array}{cc} 0&0\\
0&0\end{array}\right]}{\left[\begin{array}{cc} c&0\\
0&0\end{array}\right]})=0.
\]
Similarly, $(P_{11}+P_{21})(X)=(P_{21}+P_{22})(X)=0$, so that
$X=\left[\begin{array}{cc} 0&x\\ 0&0\end{array}\right]$.

Let $X'=\left[\begin{array}{cc} 0&0\\ 0&x\end{array}\right]$.
We claim that for any $Y\in M_2(A)$, $\tp{X}{X'}{Y}=0$.
Indeed, with $A=\left[\begin{array}{cc} a&0\\ 0&0\end{array}\right]$,
$B=\left[\begin{array}{cc} 0&b\\ 0&0\end{array}\right]$,
$C=\left[\begin{array}{cc} c&0\\ 0&0\end{array}\right]$,
we have
$A\perp X'$, $C\perp X'$ and by (\ref{eq:JP16}),
\[
\tp{X}{X'}{Y}=\tp{\tp{A}{B}{C}}{X'}{Y}=
\tp{C}{\tp{B}{A}{X'}}{Y}+\tp{A}{\tp{X'}{C}{B}}{Y}-\tp{\tp{C}{X'}{A}}{B}{Y}=0.
\]
Thus $D(X,X')=0$, which, by \cite[Lemma 1.3(a)]{FriRus85a}, 
implies that $X$ and $X'$ are orthogonal in
the Banach space sense: $\|X\pm X'\|=\max(\|X\|,\|X'\|)$. Since $
\|X+X'\|=
\|\left[\begin{array}{cc} 0&x\\ 0&x\end{array}\right]\|=\sqrt{2}\|x\|$,
 it follows that $x=0$.
The second assertion is proved similarly, using $X=
\left[\begin{array}{cc} 0&0\\ x&0\end{array}\right]$,
$X'=\left[\begin{array}{cc} 0&0\\ 0&x\end{array}\right]$.
\qed

\medskip

By interchanging rows and columns, it follows that the following triple
products all vanish (the last three by orthogonality):
\begin{equation}\label{eq:926}
\tp{\left[\begin{array}{cc} 0&0\\
0&a\end{array}\right]}{\left[\begin{array}{cc} 0&0\\
b&0\end{array}\right]}{\left[\begin{array}{cc} 0&0\\
0&c\end{array}\right]},\quad
\tp{\left[\begin{array}{cc} 0&a\\
0&0\end{array}\right]}{\left[\begin{array}{cc} b&0\\
0&0\end{array}\right]}{\left[\begin{array}{cc} 0&c\\
0&0\end{array}\right]},
\end{equation}
\begin{equation}\label{eq:927}
\tp{\left[\begin{array}{cc} 0&0\\
a&0\end{array}\right]}{\left[\begin{array}{cc} 0&0\\
0&b\end{array}\right]}{\left[\begin{array}{cc} 0&0\\
c&0\end{array}\right]},\quad
\tp{\left[\begin{array}{cc} 0&0\\
a&0\end{array}\right]}{\left[\begin{array}{cc} 0&b\\
0&0\end{array}\right]}{\left[\begin{array}{cc} 0&0\\
0&c\end{array}\right]},
\end{equation}
\begin{equation}\label{eq:928}
\tp{\left[\begin{array}{cc} 0&a\\
0&0\end{array}\right]}{\left[\begin{array}{cc} 0&0\\
b&0\end{array}\right]}{\left[\begin{array}{cc} c&0\\
0&0\end{array}\right]},\quad
\tp{\left[\begin{array}{cc} 0&0\\
0&a\end{array}\right]}{\left[\begin{array}{cc} b&0\\
0&0\end{array}\right]}{\left[\begin{array}{cc} 0&0\\
c&0\end{array}\right]}.
\end{equation}

For use in Lemma~\ref{lem:1022}, we adjoin
\[
\tp{\left[\begin{array}{cc} 0&0\\
a&0\end{array}\right]}{\left[\begin{array}{cc} b&0\\
0&0\end{array}\right]}{\left[\begin{array}{cc} 0&0\\
c&0\end{array}\right]}=
\tp{\left[\begin{array}{cc} 0&a\\
0&0\end{array}\right]}{\left[\begin{array}{cc} 0&0\\
0&b\end{array}\right]}{\left[\begin{array}{cc} 0&c\\
0&0\end{array}\right]}=0,
\]
and
\[
\tp{\left[\begin{array}{cc} 0&0\\
0&a\end{array}\right]}{\left[\begin{array}{cc} 0&b\\
0&0\end{array}\right]}{\left[\begin{array}{cc} 0&0\\
0&c\end{array}\right]}=0.
\]

\begin{lem}\label{lem:2.7} For $a,b,c\in A$, there
exists $z\in A$ such that
\[   
\tp{\left[\begin{array}{cc} a&0\\
0&0\end{array}\right]}{\left[\begin{array}{cc} 0&b\\
0&0\end{array}\right]}{\left[\begin{array}{cc} 0&c\\
0&0\end{array}\right]}
=\left[\begin{array}{cc} z&0\\ 0&0\end{array}\right].
\]  
\end{lem}
\noindent\pf\
Let $X$ denote $\tp{\left[\begin{array}{cc} a&0\\
0&0\end{array}\right]}{\left[\begin{array}{cc} 0&b\\
0&0\end{array}\right]}{\left[\begin{array}{cc} 0&c\\
0&0\end{array}\right]}$.
By  (\ref{eq:CE}), $(P_{12}+P_{22})(X)=0$ and
$(P_{12}+P_{21})(X)=0$.\qed

\begin{lem}\label{lem:2.8} For $a,b,c\in A$,
\[  
\tp{\left[\begin{array}{cc} 0&a\\
0&0\end{array}\right]}{\left[\begin{array}{cc} 0&b\\
0&0\end{array}\right]}{\left[\begin{array}{cc} 0&c\\
0&0\end{array}\right]}
=\left[\begin{array}{cc} 0&\tp{a}{b}{c}\\ 0&0\end{array}\right].
\] 
\end{lem}
\noindent\pf\
Since $P_{11}+P_{12}$ and $P_{12}+P_{22}$ are bicontractive, the
intersection of their ranges is a subtriple. Since $A$
is a \jbst\ under the product induced by $P_{12}$, and triple products are
unique, the result follows.\qed

\medskip

As noted in the proof of Lemma~\ref{lem:orth},
interchanging rows or columns 
is an isometry, hence an isomorphism. Therefore we also have, for example,
\[
\tp{\left[\begin{array}{cc} a&0\\
0&0\end{array}\right]}{\left[\begin{array}{cc} b&0\\
0&0\end{array}\right]}{\left[\begin{array}{cc} c&0\\
0&0\end{array}\right]}
=\left[\begin{array}{cc} \tp{a}{b}{c}&0\\ 0&0\end{array}\right],
\]
and so forth.

\begin{defn}\label{defn:2.7}
Define a ternary product
$\terc{a}{b}{c}$ or $\ter{a}{b}{c}$ on $A$ by
\[
\terc{a}{b}{c}=2p_{11}(
\tp{\left[\begin{array}{cc} 0&a\\
0&0\end{array}\right]}{\left[\begin{array}{cc} 0&b\\
0&0\end{array}\right]}{\left[\begin{array}{cc} c&0\\
0&0\end{array}\right]}),
\]
where $p_{11}(\left[\begin{array}{cc} a&b\\ c&d\end{array}\right])=a$.
Similarly, define two more 
ternary products $\tep{a}{b}{c}$ and $\tea{a}{b}{c}$ as follows:
\begin{equation}\label{eq:11prebis}
\tep{a}{b}{c}=2p_{11}(
\tp{\left[\begin{array}{cc} a&0\\
0&0\end{array}\right]}{\left[\begin{array}{cc} 0&0\\
b&0\end{array}\right]}{\left[\begin{array}{cc} 0&0\\
c&0\end{array}\right]})
\end{equation}
and
\begin{equation}\label{eq:11prebiss}
\tea{a}{b}{c}=2p_{11}(
\tp{\left[\begin{array}{cc} 0&0\\
c&0\end{array}\right]}{\left[\begin{array}{cc} 0&0\\
0&b\end{array}\right]}{\left[\begin{array}{cc} 0&a\\
0&0\end{array}\right]}).
\end{equation}
\end{defn}

We treat first the ternary product 
$\terc{a}{b}{c}$.
Note that, by Lemma~\ref{lem:2.7},
\begin{equation}\label{eq:11bis}
\frac{1}{2}\left[\begin{array}{cc} \terc{a}{b}{c}&0\\
0&0\end{array}\right]=\tp{\left[\begin{array}{cc} 0&a\\ 0&0
\end{array}\right]}{\left[\begin{array}{cc} 0&b\\
0&0\end{array}\right]}{\left[\begin{array}{cc} c&0\\
0&0\end{array}\right]},
\end{equation}
and that by interchanging suitable rows and columns,
\begin{eqnarray}\nonumber
\terc{a}{b}{c}&=&
2p_{21}( 
\tp{\left[\begin{array}{cc} 0&0\\
0&a\end{array}\right]}{\left[\begin{array}{cc} 0&0\\
0&b\end{array}\right]}{\left[\begin{array}{cc} 0&0\\
c&0\end{array}\right]})\\
&=&\nonumber
2p_{12}(
\tp{\left[\begin{array}{cc} a&0\\
0&0\end{array}\right]}{\left[\begin{array}{cc} b&0\\
0&0\end{array}\right]}{\left[\begin{array}{cc} 0&c\\
0&0\end{array}\right]})\\
&=&\nonumber 
2p_{22}(
\tp{\left[\begin{array}{cc} 0&0\\
a&0\end{array}\right]}{\left[\begin{array}{cc} 0&0\\
b&0\end{array}\right]}{\left[\begin{array}{cc} 0&0\\
0&c\end{array}\right]}).
\end{eqnarray}

\begin{lem}\label{1stprod} For $a,b,c\in A$,
\[ 
\terc{a}{b}{c}+\terc{c}{b}{a} = 2\tp{a}{b}{c},
\] 
and hence
\[ 
\|\terc{a}{a}{a}\|=\|a\|^3.
\]  
\end{lem}
\noindent\pf\  Given $a,b,c\in A$, it
 follows from Lemma~\ref{lem:2.3}, Lemmas~\ref{lem:2.6}-\ref{lem:2.8},
Definition~\ref{defn:2.7}
and (\ref{eq:*}) that
there are elements $x,y,z,w\in A$ such that $x+y=2\tp{a}{b}{c}$ and
\[   
\left[\begin{array}{cc} x&y\\ z&w\end{array}\right]=
\left[\begin{array}{cc} 0&0\\ 0&0\end{array}\right]
+\left[\begin{array}{cc} \ter{a}{b}{c}/2&0\\ 0&0\end{array}\right]
+\left[\begin{array}{cc} 0&\tp{a}{b}{c}\\ 0&0\end{array}\right]
+\left[\begin{array}{cc} \ter{c}{b}{a}/2&0\\ 0&0\end{array}\right].
\]  
Hence $\ter{a}{b}{c}/2+\ter{c}{b}{a}/2=x=y=\tp{a}{b}{c}$ (and $z=w=0$). 
\qed

\medskip

 We shall see later in Proposition~\ref{prodeq2} that in fact
$\ter{a}{b}{c}=\tep{a}{b}{c}=\tea{a}{b}{c}$. First we
shall show the analog of Lemma \ref{1stprod} for each of the ternary
products $\tep{a}{b}{c}$ and $\tea{a}{b}{c}$.
We note that, as above,
\begin{equation}\label{eq:14bis} 
\left[\begin{array}{cc} \tep{a}{b}{c}/2&0\\
0&0\end{array}\right]=\tp{\left[\begin{array}{cc} a&0\\
0&0\end{array}\right]}{\left[\begin{array}{cc} 0&0\\
b&0\end{array}\right]}{\left[\begin{array}{cc} 0&0\\
c&0\end{array}\right]}
\end{equation}
and
\[
\left[\begin{array}{cc} \tea{a}{b}{c}/2&0\\
0&0\end{array}\right]=\tp{\left[\begin{array}{cc} 0&0\\
c&0\end{array}\right]}{\left[\begin{array}{cc} 0&0\\
0&b\end{array}\right]}{\left[\begin{array}{cc} 0&a\\
0&0\end{array}\right]}.
\]
Moreover, by interchanging rows and/or columns,
\begin{eqnarray}
\tep{a}{b}{c}&=& \nonumber 2p_{22}( 
\tp{\left[\begin{array}{cc} 0&0\\
0&a\end{array}\right]}{\left[\begin{array}{cc} 0&b\\
0&0\end{array}\right]}{\left[\begin{array}{cc} 0&c\\
0&0\end{array}\right]})\\ 
&=&\nonumber 2p_{12}(
\tp{\left[\begin{array}{cc} 0&a\\
0&0\end{array}\right]}{\left[\begin{array}{cc} 0&0\\
0&b\end{array}\right]}{\left[\begin{array}{cc} 0&0\\
0&c\end{array}\right]})\\
&=&\nonumber 
2p_{21}(
\tp{\left[\begin{array}{cc} 0&0\\
a&0\end{array}\right]}{\left[\begin{array}{cc} b&0\\
0&0\end{array}\right]}{\left[\begin{array}{cc} c&0\\
0&0\end{array}\right]})
\end{eqnarray}
and
\begin{eqnarray}
\tea{a}{b}{c}&=&\label{eq:14prelastt} 2p_{22}(
\tp{\left[\begin{array}{cc} 0&c\\
0&0\end{array}\right]}{\left[\begin{array}{cc} b&0\\
0&0\end{array}\right]}{\left[\begin{array}{cc} 0&0\\
a&0\end{array}\right]})\\
&=&\nonumber
2p_{12}(
\tp{\left[\begin{array}{cc} 0&0\\
0&c\end{array}\right]}{\left[\begin{array}{cc} 0&0\\
b&0\end{array}\right]}{\left[\begin{array}{cc} a&0\\
0&0\end{array}\right]})\\
&=&\nonumber
2p_{21}(
\tp{\left[\begin{array}{cc}c&0\\
0&0\end{array}\right]}{\left[\begin{array}{cc}0&b\\
0&0\end{array}\right]}{\left[\begin{array}{cc} 0&0\\
0&a\end{array}\right]}).
\end{eqnarray}

\begin{prop}\label{prop:2.9}
If $A$ is an operator space such that
$M_{2}(A)_{0}$ is a bounded symmetric domain (and
consequently $M_2(A)$ and $A$ are JB*-triples), then
$\tea{a}{b}{c}+\tea{c}{b}{a}=2\tp{a}{b}{c}_A$ and
$\tep{a}{b}{c}+\tep{c}{b}{a}=2\tp{a}{b}{c}_A$.
\end{prop}
\noindent\pf\
The proof for $\tepc{\cdot}{\cdot}{\cdot}$ is similar to the proof for
$\terc{\cdot}{\cdot}{\cdot}$, using instead the identity
\begin{eqnarray*}
\tp{\left[\begin{array}{cc} a&0\\
a&0\end{array}\right]}{\left[\begin{array}{cc} 0&0\\
b&0\end{array}\right]}{\left[\begin{array}{cc}
c&0\\ c&0\end{array}\right]}&=&
\tp{\left[\begin{array}{cc} 0&0\\
c&0\end{array}\right]}{\left[\begin{array}{cc} 0&0\\
b&0\end{array}\right]}{\left[\begin{array}{cc}
0&0\\ a&0\end{array}\right]}\\
&+&
\tp{\left[\begin{array}{cc} c&0\\
0&0\end{array}\right]}{\left[\begin{array}{cc} 0&0\\
b&0\end{array}\right]}{\left[\begin{array}{cc}
0&0\\ a&0\end{array}\right]}\\
&+&
\tp{\left[\begin{array}{cc} a&0\\
0&0\end{array}\right]}{\left[\begin{array}{cc} 0&0\\
b&0\end{array}\right]}{\left[\begin{array}{cc}
c&0\\ 0&0\end{array}\right]}\\
&+&
\tp{\left[\begin{array}{cc} a&0\\
0&0\end{array}\right]}{\left[\begin{array}{cc} 0&0\\
b&0\end{array}\right]}{\left[\begin{array}{cc}
0&0\\ c&0\end{array}\right]}
\end{eqnarray*}
and the projection
\[
P(\left[\begin{array}{cc} a&b\\
c&d\end{array}\right])=\frac{1}{2}\left[\begin{array}{cc} a+c&0\\
a+c&0\end{array}\right].
\]

To prove the statement for $\teac{\cdot}{\cdot}{\cdot}$ consider (cf.
Remark~\ref{rem:2.1}) the space
\[
\t{A}=\{\t{a}=\left[\begin{array}{cc} a&a\\ a&a\end{array}\right]:a\in
A\},
\]
which is a subtriple of $M_2(A)$ since it is the range of a product $QP$
of the bicontractive projections $Q,P$ of Proposition~\ref{prop:1.3}(c). 
It follows as in
the proof of Remark~\ref{rem:2.1} that
 $\t{A}$ is a \jbst\ under the triple product $\tprime{\cdot}{\cdot}{\cdot}$
defined by
$\tprime{\t{a}}{\t{b}}{\t{c}}=4(\tp{a}{b}{c})^{\tilde{}}$.
To see this, let $D'(\t{x})\t{a}=\tprime{\t{x}}{\t{x}}{\t{a}}$ and
note that $\|\t{x}\|=2\|x\|$,
$D'(\t{x})\t{a}=4(D(x)a)^{\tilde{}}$, $e^{itD'(\t{x})}\t{y}=(
e^{4itD(x)}y)^{\tilde{}}$ and that
$(\lambda-D(\t{x}))^{-1}\t{y}=((\lambda-D(x))^{-1}y)^{\tilde{}}$.

By the uniqueness of the triple product on $M_2(A)$,
$\tp{\t{a}}{\t{b}}{\t{c}}=\tprime{\t{a}}{\t{b}}{\t{c}}$. 
Hence, by expanding
$\tp{\t{x}}{\t{y}}{\t{z}}=\tp{\left[\begin{array}{cc} x&x\\
x&x\end{array}\right]}{\left[\begin{array}{cc} y&y\\
y&y\end{array}\right]}{\left[\begin{array}{cc}
z&z\\ z&z\end{array}\right]}$ into computable terms,
\begin{eqnarray*}
\lefteqn{4\tp{x}{y}{z}^{\tilde{}}=\tp{\t{x}}{\t{y}}{\t{z}}}\\  
&=&
(\tp{x}{y}{z}+\tep{x}{y}{z}/2+\tep{z}{y}{x}/2+\ter{x}{y}{z}/2+\ter{z}{y}{x}/2
+\tea{x}{y}{z}/2+\tea{z}{y}{x}/2)^{\tilde{}}\\
&=&(3\tp{x}{y}{z}+\tea{x}{y}{z}/2+\tea{z}{y}{x}/2)^{\tilde{}}.
\end{eqnarray*}
This proves the statement for $\teac{\cdot}{\cdot}{\cdot}$. \qed

\section{Equality of the ternary products}\label{sect:equality}
In this section, we continue to assume that
 $A\subset B(H)$ is an operator space such that the open unit
ball $M_2(A)_0$ is a bounded symmetric domain.  We shall
prove the equality of the three ternary products
defined in section~\ref{sect:additivity}. Even though they agree, all
three products are needed in the proof of the crucial 
Proposition~\ref{prop:crucial}. 

In the following we shall let $\bfa\in M_2(A)$ denote 
$\left[\begin{array}{cc} a&0\\
0&a\end{array}\right]$ and $\bfab\in M_2(A)$
 denote $\left[\begin{array}{cc} 0&a\\
a&0\end{array}\right]$. By Lemmas ~\ref{lem:orth} and ~\ref{lem:2.8},
the ranges of $P_{12}+P_{21}$ and
$P_{11}+P_{22}$ are invariant under the continuous
functional calculus in a \jbst. In particular, for any $\lambda>0$,
$$\bfa^\lambda=
\left[\begin{array}{cc} a^\lambda&0\\ 0&a^\lambda\end{array}\right]
\mbox{ and }
(\bfab)^\lambda=
\left[\begin{array}{cc} 0&a^\lambda\\ a^\lambda&0\end{array}\right].
$$
Here, $\bfa^\lambda$ is defined by the triple functional calculus in
the JB*-triple $M_2(A)$ and $a^\lambda$ is defined by the triple 
functional calculus in
the JB*-triple $A$.
\begin{lem}\label{lem:new}
Let $\lambda,\mu,\nu$ be positive numbers and let $a\in A$. Then
\[
\bfa^{\lambda+\mu+\nu}=\tp{\bfa^\lambda}{\bfa^\mu}{\bfa^\nu}=
\tp{\bfab^\lambda}{\bfa^\mu}{\bfab^\nu}=
\tp{\bfab^\lambda}{\bfab^\mu}{\bfa^\nu}
\]
and
\[
\bfab^{\lambda+\mu+\nu}=\tp{\bfab^\lambda}{\bfab^\mu}{\bfab^\nu}=
\tp{\bfa^\lambda}{\bfab^\mu}{\bfa^\nu}=
\tp{\bfab^\lambda}{\bfa^\mu}{\bfa^\nu}. 
\]
\end{lem}
\noindent\pf\  
$\bfa^{\lambda+\mu+\nu}=\tp{\bfa^\lambda}{\bfa^\mu}{\bfa^\nu}$ is immediate
from the functional calculus. The proofs of the other statements are all 
proved in the same way, for example,
\begin{eqnarray*}
\lefteqn{\tp{\bfab^\lambda}{\bfa^\mu}{\bfab^\nu}=\tp
{\left[\begin{array}{cc} 0&a^{\lambda}\\ a^{\lambda}&0\end{array}\right]}
{\left[\begin{array}{cc} a^{\mu}&0\\ 0&a^{\mu}\end{array}\right]}
{\left[\begin{array}{cc} 0&a^{\nu}\\ a^{\nu}&0\end{array}\right]}=}\\
&&
\tp
{\left[\begin{array}{cc} 0&a^{\lambda}\\ 0&0\end{array}\right]}
{\left[\begin{array}{cc} a^{\mu}&0\\ 0&0\end{array}\right]}
{\left[\begin{array}{cc} 0&a^{\nu}\\ a^{\nu}&0\end{array}\right]}\\
&+&
\tp
{\left[\begin{array}{cc} 0&a^{\lambda}\\ 0&0\end{array}\right]}
{\left[\begin{array}{cc} 0&0\\ 0&a^{\mu}\end{array}\right]}
{\left[\begin{array}{cc} 0&a^{\nu}\\ a^{\nu}&0\end{array}\right]}\\
&+&
\tp
{\left[\begin{array}{cc} 0&0\\ a^{\lambda}&0\end{array}\right]}
{\left[\begin{array}{cc} a^{\mu}&0\\ 0&0\end{array}\right]}
{\left[\begin{array}{cc} 0&a^{\nu}\\ a^{\nu}&0\end{array}\right]}
\\ &+&
\tp 
{\left[\begin{array}{cc} 0&0\\ a^{\lambda}&0\end{array}\right]}
{\left[\begin{array}{cc} 0&0\\ 0&a^{\mu}\end{array}\right]}
{\left[\begin{array}{cc} 0&a^{\nu}\\ a^{\nu}&0\end{array}\right]},
\end{eqnarray*}
which further expands, using (\ref{eq:7bis})-(\ref{eq:928}) into

\begin{eqnarray*}&&
\tp
{\left[\begin{array}{cc} 0&a^{\lambda}\\ 0&0\end{array}\right]}
{\left[\begin{array}{cc} a^{\mu}&0\\ 0&0\end{array}\right]}
{\left[\begin{array}{cc} 0&0\\ a^{\nu}&0\end{array}\right]}\\
&+&
\tp
{\left[\begin{array}{cc} 0&a^{\lambda}\\ 0&0\end{array}\right]}
{\left[\begin{array}{cc} 0&0\\ 0&a^{\mu}\end{array}\right]}
{\left[\begin{array}{cc} 0&0\\ a^{\nu}&0\end{array}\right]}\\
&+&
\tp
{\left[\begin{array}{cc} 0&0\\ a^{\lambda}&0\end{array}\right]}
{\left[\begin{array}{cc} a^{\mu}&0\\ 0&0\end{array}\right]}
{\left[\begin{array}{cc} 0&a^{\nu}\\ 0&0\end{array}\right]}
\\ &+&
\tp 
{\left[\begin{array}{cc} 0&0\\ a^{\lambda}&0\end{array}\right]}
{\left[\begin{array}{cc} 0&0\\ 0&a^{\mu}\end{array}\right]}
{\left[\begin{array}{cc} 0&a^{\nu}\\ 0&0\end{array}\right]}\\
&=&
\left[\begin{array}{cc} 0&0\\
0&\tea{a^{\nu}}{a^{\mu}}{a^{\lambda}}/2\end{array}\right]
+
\left[\begin{array}{cc} \tea{a^{\lambda}}{a^{\mu}}{a^{\nu}}/2&0\\
0&0\end{array}\right]\\
&+&
\left[\begin{array}{cc} 0&0\\
0&\tea{a^{\lambda}}{a^{\mu}}{a^{\nu}}/2\end{array}\right]
+
\left[\begin{array}{cc} \tea{a^{\nu}}{a^{\mu}}{a^{\lambda}}/2&0\\
0&0\end{array}\right]
=\bfa^{\lambda+\mu+\nu}.\qed
\end{eqnarray*}   

\begin{lem}\label{Dlemma1}
$D(\bfa,\bfa)=D(\bfab,\bfab)$.
\end{lem}
\noindent\pf\
We shall use (\ref{eq:JP8}) with $z=\tp{x}{x}{y}$, which states that
\begin{equation}\label{eq:16bis}
D(\tp{x}{y}{x},\tp{x}{x}{y})=
2D(x,\tp{y}{x}{\tp{x}{x}{y}})-D(\tp{x}{\tp{x}{x}{y}}{x},y).
\end{equation}

We have, by (\ref{eq:16bis}) and Lemma~\ref{lem:new},
\begin{eqnarray*}
D(\bfa,\bfa)&=&D(\tpc{\bfab^{1/3}}{\bfa^{1/3}}{\bfab^{1/3}},\tpc{\bfab^{1/3}}{\bfab^{1/3}}{\bfa^{1/3}}\\
&=&2D(\bfab^{1/3},\tpc{\bfa^{1/3}}{\bfab^{1/3}}{\tpc{\bfab^{1/3}}{\bfab^{1/3}}{\bfa^{1/3}}})\\
&-&
D(\tpc{\bfab^{1/3}}{\tpc{\bfab^{1/3}}{\bfab^{1/3}}{\bfa^{1/3}}}{\bfab^{1/3}},\bfa^{1/3})\\
&=&2D(\bfab^{1/3},\tpc{\bfa^{1/3}}{\bfab^{1/3}}{\bfa})-
D(\tpc{\bfab^{1/3}}{\bfa}{\bfab^{1/3}},\bfa^{1/3})\\
&=&2D(\bfab^{1/3},\bfab^{5/3})-D(\bfa^{5/3},\bfa^{1/3})\\
&=&2D(\bfab,\bfab)-D(\bfa,\bfa),
\end{eqnarray*}
which proves the lemma.\qed 

\begin{lem}\label{Dlemma2}
$D(\bfa,\bfab)=D(\bfab,\bfa)$.
\end{lem}
\noindent\pf\
By Lemma~\ref{lem:new} and two applications of (\ref{eq:main}),
\begin{eqnarray*}
D(\bfa,\bfab)&=& D(\tpc{\bfa^{1/4}}{\bfab^{1/4}}{\bfab^{1/2}},\bfab)\\
&=&D(\bfab^{1/2},\{ \bfab \,\ \bfa^{1/4} \,\ \bfab^{1/4}
\})+[D(\bfa^{1/4},\bfab^{1/4}),D(\bfab^{1/2},\bfab)]\\
&=&D(\bfab^{1/2},\bfa^{3/2})+[D(\bfa^{1/4},\bfab^{1/4}),D(\bfa^{1/2},\bfa)]\\
&& (\mbox{ by Lemma~\ref{Dlemma1} since }D(\bfab^{1/2},\bfab)=D(\bfab^{3/4},\bfab^{3/4}))\\
&=&D(\bfab^{1/2},\bfa^{3/2})+D(\tp{\bfa^{1/4}}{\bfab^{1/4}}{\bfa^{1/2}},\bfa)
-D(\bfa^{1/2},\tp{\bfa}{\bfa^{1/4}}{\bfab^{1/4}})\\
&=&D(\bfab^{1/2},\bfa^{3/2})+D(\bfab,\bfa)-D(\bfa^{1/2},\bfab^{3/2}).
\end{eqnarray*}
Hence
$D(\bfa,\bfab)-D(\bfab,\bfa)=D(\bfab^{1/2},\bfa^{3/2})-D(\bfa^{1/2},\bfab^{3/2})$.

It remains to show that
$D(\bfab,\bfa^3)-D(\bfa,\bfab^3)=0$ for every $a\in  A$.

Now by (\ref{eq:JP8}) and Lemma~\ref{lem:new},
\begin{eqnarray*}
D(\bfab,\bfa^3)&=&D(\bfab,\tpc{\bfa}{\bfab}{\bfab})\\
&=& D(\tp{\bfab}{\bfa}{\bfab},\bfab)/2+D(\bfab^3,\bfa)/2\\
&=&D(\bfa^3,\bfab)/2+D(\bfab^3,\bfa)/2\\
&=&D(\bfa,\bfab^3)\mbox{ (by interchanging $\bfa$ and $\bfab$)}.
\end{eqnarray*}

This proves the lemma.\qed 

\medskip

By linearization from the preceding two lemmas we obtain
\begin{lem}\label{prodeq1}
$D(\bfa,\bfb)=D(\bfab,\bfbb)$; $D(\bfab,\bfb)=D(\bfa,\bfbb)$
\end{lem}
\noindent\pf\
From $D(\bfa+\bfb,\bfab+\bfbb)=D(\bfab+\bfbb,\bfa+\bfb)$ follows
$D(\bfb,\bfab)+D(\bfa,\bfbb)=D(\bfab,\bfb)+D(\bfbb,\bfa)$. Now replace
$a$ by $ia$  and add to obtain $D(\bfa,\bfbb)=D(\bfab,\bfb)$. The second
statement follows similarly from
$D(\bfa+\bfb,\bfa+\bfb)=D(\bfab+\bfbb,\bfab+\bfbb)$. \qed

\begin{prop}\label{prodeq2}
If $A$ is an operator space such that
$M_{2}(A)_{0}$ is a bounded symmetric domain, then
$\ter{a}{b}{c}=\tep{a}{b}{c}=\tea{a}{b}{c}$
\end{prop}
\noindent\pf\
By expanding as in the second part of the proof of Lemma~\ref{Dlemma2},  
\begin{eqnarray*}
D(\bfa,\bfb)
\left[\begin{array}{cc} x&0\\ 0&0\end{array}\right]&=&
\tp
{\left[\begin{array}{cc} a&0\\ 0&a\end{array}\right]}
{\left[\begin{array}{cc} b&0\\ 0&b\end{array}\right]}
{\left[\begin{array}{cc} x&0\\ 0&0\end{array}\right]}\\
&=&
\tp
{\left[\begin{array}{cc} a&0\\ 0&0\end{array}\right]}
{\left[\begin{array}{cc} b&0\\ 0&0\end{array}\right]}
{\left[\begin{array}{cc} x&0\\ 0&0\end{array}\right]}\\
&=&\left[\begin{array}{cc} \tp{a}{b}{x}&0\\ 0&0\end{array}\right]
\end{eqnarray*}
and
\begin{eqnarray*}
D(\bfab,\bfbb)
\left[\begin{array}{cc} x&0\\ 0&0\end{array}\right]
&=&
\tp
{\left[\begin{array}{cc} 0&a\\ a&0\end{array}\right]}
{\left[\begin{array}{cc} 0&b\\ b&0\end{array}\right]}
{\left[\begin{array}{cc} x&0\\ 0&0\end{array}\right]}\\
&=&
\tp
{\left[\begin{array}{cc} 0&a\\ 0&0\end{array}\right]}
{\left[\begin{array}{cc} 0&b\\ 0&0\end{array}\right]}
{\left[\begin{array}{cc} x&0\\ 0&0\end{array}\right]}\\
&+&
\tp
{\left[\begin{array}{cc} 0&0\\ a&0\end{array}\right]}
{\left[\begin{array}{cc} 0&0\\ b&0\end{array}\right]}
{\left[\begin{array}{cc} x&0\\ 0&0\end{array}\right]}\\
&=&
\left[\begin{array}{cc} \ter{a}{b}{x}/2+\tep{x}{b}{a}/2&0\\
0&0\end{array}\right],
\end{eqnarray*} so that $\ter{x}{b}{a}=\tep{x}{b}{a}$.

Similarly,
\begin{eqnarray*}
\left[\begin{array}{cc} 0&\tea{x}{b}{a}/2\\
\tea{a}{b}{x}/2&0\end{array}\right]&=&
\tp
{\left[\begin{array}{cc} a&0\\ 0&a\end{array}\right]}
{\left[\begin{array}{cc} 0&b\\ b&0\end{array}\right]}
{\left[\begin{array}{cc} x&0\\ 0&0\end{array}\right]}
=D(\bfa,\bfbb)
\left[\begin{array}{cc} x&0\\ 0&0\end{array}\right]
\\
&=&
D(\bfab,\bfb)
\left[\begin{array}{cc} x&0\\ 0&0\end{array}\right]
=
\tp
{\left[\begin{array}{cc} 0&a\\ a&0\end{array}\right]}
{\left[\begin{array}{cc} b&0\\ 0&b\end{array}\right]}
{\left[\begin{array}{cc} x&0\\ 0&0\end{array}\right]}\\
&=&
\left[\begin{array}{cc} 0&\ter{x}{b}{a}/2\\
\tep{a}{b}{x}/2&0\end{array}\right],
\end{eqnarray*}
so that $\tea{x}{b}{a}=\ter{x}{b}{a}$. \qed

\section{Main result}\label{sect:main}

\begin{prop}\label{prop:crucial}
Let $X$ be an operator space such that
$M_{2}(X)_{0}$ is a bounded symmetric domain. Then
$(X,\ter{\cdot}{\cdot}{\cdot},\|\cdot\|)$ is a $C^*$-ternary ring in the
sense of Zettl \cite{Zettl83} {\rm (see Definition~\ref{defn:1.2})}
and its \jbst\ product {\rm 
(see the beginning of section~\ref{sect:additivity})}
satisfies
$\tp{a}{b}{c}=(\ter{a}{b}{c}+\ter{c}{b}{a})/2$.
\end{prop}
\noindent\pf\
It was already shown in Lemma \ref{1stprod} that
$\tp{a}{b}{c}=(\ter{a}{b}{c}+\ter{c}{b}{a})/2$ and that
$\|\ter{a}{a}{a}\|=\|a\|^3$ and it is clear that
$\|\ter{a}{b}{c}\|\le\|a\|\|b\|\|c\|$.

It remains to show associativity. To prove this we will use 
Lemma~\ref{lem:orth}
and Proposition~\ref{prodeq2}.  For $a,b,c,d,e\in X$,
let
\[A=\left[\begin{array}{cc} 0&a\\
0&0\end{array}\right],B=\left[\begin{array}{cc} 0&b\\
0&0\end{array}\right],
C=\left[\begin{array}{cc} c&0\\
0&0\end{array}\right],D=\left[\begin{array}{cc} 0&0\\
d&0\end{array}\right],
E=\left[\begin{array}{cc} 0&0\\ e&0\end{array}\right].
\]
 Then
\begin{eqnarray*}
\ter{\ter{a}{b}{c}}{d}{e}&=&\tep{\ter{a}{b}{c}}{d}{e}
=2p_{11}(\tpc
{\left[\begin{array}{cc} \ter{a}{b}{c}&0\\ 0&0\end{array}\right]}
{\left[\begin{array}{cc} 0&0\\ d&0\end{array}\right]}
{\left[\begin{array}{cc} 0&0\\ e&0\end{array}\right]})
\mbox{\quad (by (\ref{eq:11prebis}))}
\\
&=&4p_{11}(\tpc{\tp
{\left[\begin{array}{cc} 0&a\\ 0&0\end{array}\right]}
{\left[\begin{array}{cc} 0&b\\ 0&0\end{array}\right]}
{\left[\begin{array}{cc} c&0\\ 0&0\end{array}\right]}}
{\left[\begin{array}{cc} 0&0\\
d&0\end{array}\right]}{\left[\begin{array}{cc} 0&0\\
e&0\end{array}\right]})
\\ &&\mbox{ (by (\ref{eq:11bis}))}
\\
&=&4p_{11}(\tp{E}{D}{\tp{C}{B}{A}})\quad\mbox{ (by commutativity of the triple
product)}\\
&=&4p_{11}(\tp{C}{B}{\tp{E}{D}{A}})+
4p_{11}(\tp{\tp{E}{D}{C}}{B}{A})\\
&&-
4p_{11}(\tp{C}{\tp{B}{E}{D}}{A})\quad\mbox{ (by (\ref{eq:JP15}))}\\
&=&
0+
4p_{11}(\tpc{\tp
{\left[\begin{array}{cc} 0&0\\ e&0\end{array}\right]}
{\left[\begin{array}{cc} 0&0\\ d&0\end{array}\right]}
{\left[\begin{array}{cc} c&0\\ 0&0\end{array}\right]}}
{\left[\begin{array}{cc} 0&b\\
0&0\end{array}\right]}{\left[\begin{array}{cc} 0&a\\
0&0\end{array}\right]})+0\\
&=&2p_{11}(  
\tp
{\left[\begin{array}{cc} 0&a\\ 0&0\end{array}\right]}
{\left[\begin{array}{cc} 0&b\\ 0&0\end{array}\right]}
{\left[\begin{array}{cc} \tep{c}{d}{e}&0\\ 0&0\end{array}\right]}
)
\mbox{ (by (\ref{eq:14bis}))}\\  
&=&\ter{a}{b}{\tep{c}{d}{e}}=\ter{a}{b}{\ter{c}{d}{e}}.
\end{eqnarray*}

To complete the proof of associativity,  consider
\begin{eqnarray*}
\ter{a}{\ter{d}{c}{b}}{e}&=&\tea{a}{\tea{d}{c}{b}}{e}\\
&=&2p_{11}(\tp
{\left[\begin{array}{cc} 0&a\\ 0&0\end{array}\right]}
{\left[\begin{array}{cc} 0&0\\ 0&\tea{d}{c}{b}\end{array}\right]}
{\left[\begin{array}{cc} 0&0\\ e&0\end{array}\right]})
\mbox{ (by (\ref{eq:11prebiss}))}\\ 
&=&4p_{11}(\tpc
{\left[\begin{array}{cc} 0&a\\ 0&0\end{array}\right]}
{\tp
{\left[\begin{array}{cc} 0&0\\ d&0\end{array}\right]}
{\left[\begin{array}{cc} c&0\\ 0&0\end{array}\right]}
{\left[\begin{array}{cc} 0&b\\ 0&0\end{array}\right]}}
{\left[\begin{array}{cc} 0&0\\ e&0\end{array}\right]})
\\ &&\mbox{ (by (\ref{eq:14prelastt}))}\\
\\
&=&4p_{11}(\tp{A}{\tp{D}{C}{B}}{E})\\
&=&4p_{11}(\tp{\tp{A}{B}{C}}{D}{E})+4p_{11}(\tp{E}{B}{\tp{A}{D}{C}})-
4p_{11}(\tp{C}{\tp{B}{A}{D}}{E})
\\ &&\mbox{ (by (\ref{eq:JP16}))}\\  
&=&4p_{11}(\tp{\tp{A}{B}{C}}{D}{E})
\mbox{ (since $A\perp D$)}\\
&=&4p_{11}(\tpc{\tp
{\left[\begin{array}{cc} 0&a\\ 0&0\end{array}\right]}
{\left[\begin{array}{cc} 0&b\\ 0&0\end{array}\right]}
{\left[\begin{array}{cc} c&0\\ 0&0\end{array}\right]}}
{\left[\begin{array}{cc} 0&0\\
d&0\end{array}\right]}{\left[\begin{array}{cc} 0&0\\
e&0\end{array}\right]})\\
&=&2p_{11}(\tp
{\left[\begin{array}{cc} \ter{a}{b}{c}&0\\ 0&0\end{array}\right]}
{\left[\begin{array}{cc} 0&0\\ d&0\end{array}\right]}
{\left[\begin{array}{cc} 0&0\\ e&0\end{array}\right]})
=\tep{\ter{a}{b}{c}}{d}{e}
\mbox{ (by (\ref{eq:14bis}))}\\
&=&\ter{\ter{a}{b}{c}}{d}{e}.\qed
\end{eqnarray*} 


\begin{lem}\label{lem:1022}
Let $A$ be an operator space such that
$M_{2}(A)_{0}$ is a bounded symmetric domain, so that
by Proposition~\ref{prop:crucial}, $A$ is a C*-ternary ring. 
Suppose that the C*-ternary
ring $A$ is isomorphic to a TRO, that is, $A_{-1}=0$ in
Theorem~\ref{thm:1}. Form the ternary product 
$\ter{\cdot}{\cdot}{\cdot}_{M_2(A)}$
induced by the ternary product on $A$ as if it was
 ordinary matrix multiplication, that is, if
$X=[x_{ij}],\ Y=[y_{kl}],\ Z=[z_{pq}]\in M_2(A)$, then
$\ter{X}{Y}{Z}_{M_2(A)}$ is the matrix whose $(i,j)$-entry
is $\sum_{p,q}\ter{x_{ip}}{y_{qp}}{z_{qj}}$.
Then
\[
2\tp{X}{Y}{Z}_{M_2(A)}=
\ter{X}{Y}{Z}_{M_2(A)}+
\ter{Z}{Y}{X}_{M_2(A)}
\]
\end{lem}
\pf\
It suffices to prove that $\tp{X}{X}{X}_{M_2(A)}=
\ter{X}{X}{X}_{M_2(A)}$. 

In the first place,
\[
\ter{X}{X}{X}_{M_2(A)}=
\left[\begin{array}{c|c}          
&\\
\ter{x_{11}}{x_{11}}{x_{11}}+
\ter{x_{12}}{x_{12}}{x_{11}}
&
\ter{x_{11}}{x_{11}}{x_{12}}+
\ter{x_{12}}{x_{12}}{x_{12}}\\
&\\
+\ter{x_{11}}{x_{21}}{x_{21}}+
\ter{x_{12}}{x_{22}}{x_{21}}
&
+\ter{x_{11}}{x_{21}}{x_{22}}+
\ter{x_{12}}{x_{22}}{x_{22}}\\
&\\\hline
&\\
\ter{x_{21}}{x_{11}}{x_{11}}+
\ter{x_{22}}{x_{12}}{x_{11}}
&
\ter{x_{21}}{x_{11}}{x_{12}}+
\ter{x_{22}}{x_{12}}{x_{12}}\\
&\\
+\ter{x_{21}}{x_{21}}{x_{21}}+
\ter{x_{22}}{x_{22}}{x_{21}}
&
+\ter{x_{21}}{x_{21}}{x_{22}}
+\ter{x_{22}}{x_{22}}{x_{22}}\\
&
\end{array}\right]
\]

On the other hand, by using 
Lemmas~\ref{lem:orth},\ref{lem:2.6},\ref{lem:2.8} and \ref{1stprod}, and
Proposition~\ref{prop:2.9},
\[
\sum_{k,l,p,q}\tp{P_{11}(X)}{P_{kl}(X)}{P_{pq}(X)}=
\left[\begin{array}{c|c}          
&\\
\tp{x_{11}}{x_{11}}{x_{11}}+
\ter{x_{12}}{x_{12}}{x_{11}}/2
&
\ter{x_{11}}{x_{11}}{x_{12}}/2\\
&\\
+\ter{x_{11}}{x_{21}}{x_{21}}/2
&
+\ter{x_{11}}{x_{21}}{x_{22}}/2\\
&\\\hline
&\\
\ter{x_{21}}{x_{11}}{x_{11}}/2+
\ter{x_{22}}{x_{12}}{x_{11}}/2
&
0
\\
&
\end{array}\right],
\]
\[
\sum_{k,l,p,q}\tp{P_{12}(X)}{P_{kl}(X)}{P_{pq}(X)}=
\left[\begin{array}{c|c}          
&\\
\ter{x_{12}}{x_{22}}{x_{21}}/2
&
\tp{x_{12}}{x_{12}}{x_{12}}+\ter{x_{11}}{x_{11}}{x_{12}}/2\\
&\\
+\ter{x_{12}}{x_{12}}{x_{11}}/2
&
+\ter{x_{12}}{x_{22}}{x_{22}}/2\\
&\\\hline
&\\
0&
\ter{x_{21}}{x_{11}}{x_{12}}/2
+
\ter{x_{22}}{x_{12}}{x_{12}}/2
\\
&  
\end{array}\right],
\]
\[
\sum_{k,l,p,q}\tp{P_{21}(X)}{P_{kl}(X)}{P_{pq}(X)}=
\left[\begin{array}{c|c}          
&\\
\ter{x_{11}}{x_{21}}{x_{21}}/2+
\ter{x_{12}}{x_{22}}{x_{21}}/2&0\\
&\\\hline
&\\
\tp{x_{21}}{x_{21}}{x_{21}}+\ter{x_{21}}{x_{11}}{x_{11}}/2&
\ter{x_{21}}{x_{11}}{x_{12}}/2\\
&\\
+\ter{x_{22}}{x_{22}}{x_{21}}/2&
+
\ter{x_{21}}{x_{21}}{x_{22}}/2
\\
&  
\end{array}\right],
\]
and
\[
\sum_{k,l,p,q}\tp{P_{22}(X)}{P_{kl}(X)}{P_{pq}(X)}=
\left[\begin{array}{c|c}          
&\\
0&\ter{x_{11}}{x_{21}}{x_{22}}/2+
\ter{x_{12}}{x_{22}}{x_{22}}/2\\
&\\\hline
&\\
\ter{x_{22}}{x_{12}}{x_{11}}/2&
\tp{x_{22}}{x_{22}}{x_{22}}
+\ter{x_{21}}{x_{21}}{x_{22}}/2\\
&\\
+\ter{x_{22}}{x_{22}}{x_{21}}/2&
\ter{x_{22}}{x_{12}}{x_{12}}/2
\\
&  
\end{array}\right].
\]

Since 
$
\tp{X}{X}{X}_{M_2(A)}=\sum_{i,j}\sum_{k,l,p,q}
\tp{P_{ij}(X)}{P_{kl}(X)}{P_{pq}(X)}
$
and 
$\tp{x}{x}{x}=\ter{x}{x}{x}$, the lemma follows.
\qed

\medskip

We now state and prove the main result of this paper.

\begin{thm}\label{thm:first}
Let $A\subset B(H)$ be an operator space and suppose that $M_n(A)_0$
is a bounded symmetric domain for some 
$n \geq 2$. Then $A$ is n-isometric to a
ternary ring of
operators (TRO). If $M_n(A)_0$ is a 
bounded symmetric domain for all $n\ge 2$, then
$A$ is ternary isomorphic and completely isometric to a TRO.
\end{thm}
\noindent\pf\  
The second statement follows from the first one.
Suppose $n=2$.
From Theorem~\ref{thm:1} 
and Proposition~\ref{prop:crucial},   
we know that $A=A_{1} \oplus A_{-1}$ where
$A_{1}$ is ternary isomorphic  to a TRO $B$
and $A_{-1}$ is anti-isomorphic to a TRO $C$. 
Let $\varphi:A_{-1}\rightarrow C$ be an anti-isomorphism.
Since $C$ is a
JB*-triple under the product $\{ x \,\ y \,\ z \} =
(1/2)(xy^{\ast}z+zy^{\ast}x)$ and 
%
$\varphi$ is an isometry, hence a triple isomorphism, it follows that
\[
\varphi(x)\varphi(x)^*\varphi(x)=\varphi\tp{x}{x}{x}=\varphi\ter{x}{x}{x}=
-\varphi(x)\varphi(x)^*\varphi(x)
\]
so that $\varphi(x)\varphi(x)^*\varphi(x)=0$ and $x=0$.
Thus $A_{-1}=0$ and
$A$ is ternary isomorphic to a TRO $B$. Let $\psi:A\rightarrow
B$ be a surjective ternary isomorphism. Then by Lemma~\ref{lem:1022},
the amplification $\psi_2$
is a triple isomorphism of the JB*-triple $M_2(A)$ onto the JB*-triple
$M_2(B)$, with the triple product 
\[
\tp{R}{S}{T}_{M_2(B)}:=(RS^*T+TS^*R)/2,
\]
implying that $\psi_2$ is a triple isomorphism, hence an isometry.
Thus, $A$ is 2-isometric
to $B$, proving the theorem for $n=2$.

The general case for $M_{n}(A)$ is now not difficult to obtain.
We require only one short lemma. 

\begin{lem}\label{arm}
Let $A$ be an operator space such that for some $n\ge 3$,
 $M_{n}(A)$ has a JB*-triple
structure. Then for $X,Y,Z\in M_n(A)$,  the following products all 
vanish:
\begin{itemize}
\item $\{P_{ij}(X) \,\ P_{kj}(Y) \,\ P_{lj}(Z)\}$ (for distinct $i,k,l$)
\item
$\{P_{ij}(X) \,\ P_{ik}(Y) \,\ P_{il}(Z)\}$ (for distinct $j,k,l$)
\item
 $\{P_{ij}(X) \,\ P_{kl}(Y) \,\ P_{pq}(Z)\}$ 
(for $i\ne k$, $j\ne l$ and either $p\not\in\{i,k\}$ or $q\not \in\{j,l\}$
\end{itemize}  
\end{lem}
\noindent\pf\  
Two applications of the fact that the range of
a bicontractive projection on a JB*-triple is a subtriple yield that
$\{P_{ij}(X) \,\ P_{kj}(Y) \,\ P_{lj}(Z)\}$ lies in
$(P_{ij}+P_{kj}+P_{lj})M_{n}(A)$. However, by a conditional expectation
property,
\[
(P_{ij}+P_{kj})\{P_{ij}(X) \,\ P_{kj}(Y) \,\ P_{lj}(Z)\}
=(P_{ij}+P_{kj})\{P_{ij}(X)
\,\ P_{kj}(Y) \,\ 0\}=0.
\]
A similar calculation shows $(P_{kj}+P_{lj})\{P_{ij}(X) \,\ P_{kj}(Y) \,\
P_{lj}(Z)\}=0$, proving the first statement. A similar agrument proves
the second statement. The proof of the last statement is
the same as the proof of Lemma~\ref{lem:orth}. 
For $n=3$, one needs to prove, for example, that
\[
D(\left[\begin{array}{ccc} a &0&0\\
0&0&0\\
0&0&0
 \end{array}\right]
,
\left[\begin{array}{ccc} 0 &0&0\\
0&b&0\\
0&0&0
 \end{array}\right]
)=0 \qed
\]

\medskip

Returning to the proof of Theorem~\ref{thm:first},
if $M_n(A)$ is a JB*-triple, then $M_{2}(A)$, which is isometric to the
range of a contractive projection on $M_n(A)$, is also a JB*-triple. Hence, by
the $n=2$ case, $A$ is a C*-ternary ring which is ternary isomorphic and
isometric under a map $\phi$ to a TRO $B$ and $M_{2}(A)$
is triple isomorphic and isometric to $M_{2}(B)$ under the amplification
$\phi_{2}$. Every triple product $\tp{X}{Y}{Z}$
in $M_n(A)$ is the sum of products of 
the form 
$\{P_{ij}(X) \,\ P_{kl}(Y) \,\ P_{pq}(Z)\}$.
By Lemma \ref{arm}, every such product of
matrix elements in $M_{n}(A)$ is either zero or takes place in the
intersection of two rows with two columns. The subspace of $M_n(A)$
defined by one such  intersection is
a subtriple of $M_n(A)$ since it is the range of the product of two
bicontractive projections.  It is
 isometric, via
\[
P_{ij}(X)+P_{il}(Y)+P_{kj}(Z)+P_{kl}(W)\mapsto 
\left[\begin{array}{cc} P_{ij}(X) &P_{il}(Y)\\
P_{kj}(Z) &P_{kl}(W) \end{array}\right],
\] hence triple isomorphic,
 to $M_{2}(A)$. Hence, by the proof of the $n=2$ case, 
all triple products in $M_n(A)$ are the natural ones obtained
from the ternary structure on $A$ as in Lemma~\ref{lem:1022}.
It follows
 that $M_{n}(A)$ is triple isomorphic to $M_{n}(B)$
via the amplification map $\phi_{n}$ which is thus an isometry. \qed

\medskip

As application, we offer two corollaries.

\begin{cor}\label{cor:4.4}
Let ${\mathcal A}\subset B(H,K)$ be a TRO and let $P$
be a completely contractive projection on $\mathcal A$. Then the range of
$P$ is completely isometric to another TRO.
\end{cor}
\noindent\pf\
Since $\mathcal A$ is a TRO, $M_n({\mathcal A})$ is a JB*-triple. Therefore
$M_n(P({\mathcal A}))=P_n(M_n({\mathcal A}))$ is a JB*-triple, and its
unit ball is a bounded symmetric domain. \qed

\medskip

Another way to obtain this corollary is to note that every TRO is a corner
of a C*-algebra and hence the range of a completely contractive projection
on that algebra.  By composing these two projections, the corollary is
reduced to \cite{Youngson83}.

\medskip

Our second corollary is a variant of the fundamental Choi-Effros result.

\begin{cor}\label{cor:4.6}
Let $P$ be a unital $2$-positive projection on a unital JC*-algebra $A$.
Then $P(A)$ is $2$-isometric to a C*-algebra.  If $P$
is completely positive and unital, then 
$P(A)$ is  completely isometric to a C*-algebra. 
\end{cor}

In order to state our second theorem, we recall 
that a complex Banach space $A$ is linearly isometric to 
a unital JB*-algebra if and only if its open unit ball $A_0$ is a bounded
symmetric domain of tube type \cite{BraKauUpm78}. In \cite{Upmeier84},
a necessary and sufficient condition, involving the Lie algebra of all
complete holomorphic vector fields on $A_0$, is given for such $A$ to be
obtained from a
C*-algebra with the anticommutator product. Our next theorem gives 
a holomorphic characterization of C*-algebras up to {\it complete} isometry.

\begin{thm}\label{thm:4.5}
Let $A\subset B(H)$ be an operator space and suppose that $M_n(A)_0$
is a bounded symmetric domain for some
$n \geq 2$. If the induced bounded symmetric
domain structure on $A_0$ is of tube 
type, then $A$ is n-isometric to a
C*-algebra. If $M_n(A)_0$ is a bounded symmetric domain 
for all $n\ge 2$ and $A_0$ is of tube type, then 
$A$ is completely isometric to a C*-algebra.
\end{thm}
\noindent\pf\
By Theorem~\ref{thm:first}, we may assume that $A$ is a TRO.
Since $A$ has the
structure of a unital JB*-algebra, there is a partial isometry $u\in A$
such that $au^*u=uu^*a=a$ for every $a\in A$. Then $A$ becomes a C*-algebra
with product $a\cdot b=au^*b$ and involution $a^\sharp=ua^*u$. Since
$ab^*c=a\cdot b^\sharp\cdot c$, and ternary isomorphisms of TRO's are
complete isometries, the result follows. \qed

\begin{rem}
One can construct operator spaces that are 2-isometric to a C*-algebra $A$ which
are not completely isometric to $A$.   Hence, if $M_2(A)_0$ is a bounded symmetric domain
it does not follow that $M_n(A)_0$ is a bounded symmetric domain for every $n\ge2$. It
would be interesting to see if this were true under some further condition on $A$.

The proof of Theorem~\ref{thm:first} seems to require a bounded symmetric
domain structure
on $M_2(A)_0$, not simply on $M_{1,2}(A)_0$ for example. It would be interesting
to see what could be said if it is assumed that $M_{1,n}(A)_0$ were a bounded 
symmetric domain for every $n\ge 2$.
\end{rem}

\end{document}